\documentclass[12pt,a4paper]{article}
\usepackage[textwidth=16cm, textheight=23cm,marginratio=1:1]{geometry}
\usepackage[T1]{fontenc}            
\usepackage[utf8]{inputenc}
\usepackage{amsmath}
\usepackage{amsfonts}
\usepackage{mathtools}
\usepackage{amssymb}
\usepackage{comment}
\usepackage{tikz-cd}
\usepackage{pstricks,xy,pst-node}
\usepackage{graphicx}
\xyoption{all}
\usepackage{amsthm}
\usepackage{booktabs} 
\usepackage{mathrsfs} 
\usepackage{url}
\usepackage{caption} 
\usepackage{enumitem}
\usepackage{sidecap}

\usepackage[breaklinks,bookmarks=false, hidelinks]{hyperref} 
\hypersetup{
    colorlinks,
    linkcolor={red}, citecolor={teal}
}

\makeatletter
\let\c@table\c@figure 
\let\ftype@table\ftype@figure 
\makeatother

\newtheorem{thm}[equation]{Theorem}

\newtheorem{lem}[equation]{Lemma}
\newtheorem{conj}[equation]{Conjecture}
\newtheorem*{thm*}{Theorem}

\newtheorem*{problem*}{Problem}

\newtheorem{funfact}[equation]{Fun Fact}

\theoremstyle{remark}

\theoremstyle{definition}

\DeclareMathOperator{\socdim}{socdim}

\DeclareMathOperator{\Gr}{Gr}

\newcommand{\Hilb}{{\ccH}ilb}

\newcommand{\PP}{\mathbb{P}}

\def\GG{{\mathbb G}}

\def\PP{{\mathbb P}}

\def\kk{{\Bbbk}}
\newcommand{\cactus}[2]{\mathfrak{K}_{#1}\left( #2 \right)}
\newcommand{\cactusop}[2]{\mathfrak{K}^{\circ}_{#1}\left( #2 \right)}

\newcommand{\PPof}[1]{\PP\left(#1\right)}

\newcommand{\ccH}{{\mathcal{H}}}
\newcommand{\ccI}{{\mathcal{I}}}

\newcommand{\ccO}{{\mathcal{O}}}

\newcommand{\ccR}{{\mathcal{R}}}

\newcommand{\gotm}{\mathfrak{m}}

\newcommand{\gots}{\mathfrak{s}}

\DeclareMathOperator{\cchar}{char}
\DeclareMathOperator{\Spec}{Spec}

\DeclareMathOperator{\Supp}{Supp} 

\newcommand{\Gor}{\operatorname{Gor}}

\newcommand{\set}[1]{\left\{#1\right\}}

\newcommand{\linspan}[1]{\left\langle#1\right\rangle}

\newcommand{\reduced}[1]{{#1}_{\operatorname{red}}}

\newenvironment{prf}[1][]
  {\medskip\par\noindent{\bf Proof#1. }}
  {\nopagebreak\qed\par\smallskip}

\newcommand{\noprf}{\nopagebreak\qed\par}

\newcounter{betweenenumi}

\newcommand{\eemail}[1]{{\href{mailto:#1}{\nolinkurl{#1}}}}

\title{Grassmann cactus variety and socle dimension}

\author{Weronika Buczy{\'n}ska \and Jaros{\l}aw Buczy{\'n}ski \and Maciej Ga{\l}{\k a}zka}

\usepackage[textsize=tiny]{todonotes}

\begin{document}
\date{September 8, 2025}

\maketitle

\begin{abstract}
Grassmann cactus variety is a common generalisation of Grassmann secant variety and cactus variety.
In their definitions one considers the vector spaces of fixed dimension that are contained in the linear span of some finite schemes.
We prove that to characterise Grassmann cactus varieties it is enough to use finite schemes that locally have low socle dimension.
This motivates the study of parameter spaces of such schemes and simplifies calculations of examples of Grassmann cactus varieties.
\end{abstract}

\medskip
{\footnotesize
\noindent\textbf{addresses:} \\
W.~Buczy\'nska, \eemail{wkrych@mimuw.edu.pl},
   Faculty of Mathematics, Computer Science and Mechanics, University of Warsaw, ul. Banacha 2, 02-097 Warsaw, Poland\\
J.~Buczy\'nski, \eemail{jabuczyn@impan.pl},
   Institute of Mathematics of the Polish Academy of Sciences, ul. \'Snia\-dec\-kich 8, 00-656 Warsaw, Poland\\
M.~Ga{\l}{\k a}zka, \eemail{mgalazka@dgist.ac.kr},
  Daegu-Gyeongbuk Institute of Science and Technology, 333 Techno Jun\-gang-daero, Hyeonpung-eup, Dalseong-gun, Daegu 42988, South
   Korea, and Warsaw University of Life Sciences, ul. Nowoursynowska 166, 02-787 Warszawa, Poland

\noindent\textbf{keywords:}
Grassmann cactus variety, socle dimension, finite schemes, secant variety, linear span

\noindent\textbf{AMS Mathematical Subject Classification 2020:}
Primary: 14A15, 
Secondary:
13H10, 
14D99, 
14N07 
}


\section{Introduction}\label{sec_intro}
Throughout this note, we work over an algebraically closed field $\kk$ of any characteristic.

Suppose $V$ is a $\kk$-vector space, and let $\PP V$ be the corresponding (naive) projective space $\PP V = (V \setminus \set{0})/\GG_{\gotm}$ over $\kk$.
Let $\Gr(k, V)$ be the Grassmannian of $k$-dimensional linear subspaces of $V$, so that $\Gr(1, V) = \PP V$, and $\Gr(\dim V-1, V) = \PP V^*$.
If $R\subset \PP V$ is any closed subscheme, then by linear span $\linspan{R} \subset V$ we mean the smallest linear subspace of $V$ such that $R\subset \PP\linspan{R}$ as a scheme.
Consider a fixed quasiprojective scheme $X$ with a fixed locally closed embedding $X\hookrightarrow \PP V$.

We define the \emph{Grassmann cactus variety} of $X$ to be:
\[
 \cactus{r,k}{X} := \overline{\set{E\in \Gr(k, V) \mid
    \exists \text{ subscheme } R\subset X \text{ s.t.} \dim R = 0, \deg R \leqslant r, E \subset \linspan{R} }}.
\]
Here and throughout the article the overline $\overline{\set{\dotsc}}$ denotes the Zariski closure of the set $\set{\dotsc}$. Note that if $X$ has at least $r$ distinct points of support (for instance, if $\dim X >0$), then in the definition we can assume $\deg R = r$ instead of $\deg R \leqslant r$.
The case of $k=1$ is the well established in the literature case of \emph{cactus variety} $\cactus{r}{X}\subset \PP V$,
  see for instance \cite{nisiabu_jabu_cactus}, \cite{galazka_vb_cactus}, \cite{nisiabu_jabu_farnik_cactus_Fujita},
  \cite{bernardi_ranestad_cactus_rank_of_cubics},
  \cite{bernardi_jelisiejew_marques_ranestad_cactus_rank_of_a_general_form}.
Also if in the definition we only consider smooth schemes $R$ (that is, unions of $r$ distinct points of $X$) then we obtain the \emph{Grassmann secant varieties} \cite{chiantini_coppens_Grassmannians_of_secant_varieties}, \cite{ciliberto_cools_Grassmann_secant_extremal_varieties}, \cite{landsberg_jabu_ranks_of_tensors}.
Thus the definition of $\cactus{r,k}{X}$ is a common generalisation of cactus varieties and Grassmann secant varieties, which also appears in \cite{galazka_mandziuk_rupniewski_distinguishing},
\cite{nisiabu_jabu_abcd}.
Furthermore, we expect that appropriately exploited $\cactus{r,k}{X}$ will lead to cactus barriers on determinantal methods for border rank, which are stronger than  \cite{efremenko_garg_oliveira_wigderson_barriers_for_rank_methods} or  \cite{garg_makam_oliveira_wigderson_more_barriers_for_rank_methods}.

In this note, the term ``socle'' plays a central role. Let $A$ be a $0$-dimensional local ring with maximal ideal $\gotm$. The socle of $A$ is the annihilator of $\gotm$:
\[
  \gots = \set{a \in A \mid a\cdot \gotm  = 0}.
\]
The socle $\gots$ of $A$ is at the same time an ideal in $A$, an $A$-module (since all ideals are modules) and a $\kk$-vector space (since $\gotm$ acts trivially on $\gots$, and therefore $\gots$ is an $A/\gotm$-module, where $A/\gotm\simeq \kk$).
See \cite[Sec.~21.1]{eisenbud} for a geometric explanation of the name and a discussion of the notion of socle in a slightly more general context.

Any finite scheme $R$ over $\kk$ can be written as $R \simeq \bigsqcup_{i=1}^{\#\Supp R} \Spec A_i$ with each $A_i$ a finite local algebra over $\kk$.
We say $R$ has socle dimension at most $l$ (denoted by $\socdim(R)\leqslant l$) if for all $i$ the dimension (as a $\kk$-vector space) of the socle of $A_i$ is at most $l$.
In particular, $R$ is Gorenstein if and only if $R$ has socle dimension at most $1$ \cite[Prop.~21.5(a)\&(c)]{eisenbud}.

Our result is the following:
\begin{thm}
   \label{thm_cactus_and_socle}
   In the definition of $\cactus{r,k}{X}$ it is enough to use schemes of socle dimension at most $k$:
   \[
      \cactus{r,k}{X} = \overline{\set{E\in \Gr(k, V) \mid
      \exists R\subset X \text{ s.t.} \dim R = 0, \deg R \leqslant r, \socdim R\leqslant k, E \subset \linspan{R} }}.
   \]
\end{thm}

Case $k=1$ of this theorem is
\cite[Prop.~2.2(ii)]{nisiabu_jabu_cactus}.
This special case is often used in the context of comparing secant varieties or varieties of sums of powers (VSP) with subsets of Hilbert schemes.
The reduction to Gorenstein subschemes may reduce the complexity of the problem, as Gorenstein schemes are slightly better behaved than general schemes, or at least pathologies appear a little later.
See for instance~\cite[Rem.~3.3]{ranestad_voisin_VSP_and_divisors_in_the_moduli_of_cubic_fourfolds}, \cite[\S1.1]{jelisiejew_ranestad_schreyer_VPS2},
\cite[\S1]{bernardi_ranestad_cactus_rank_of_cubics},
\cite[Thm~1.2]{jabu_keneshlou_cactus_scheme},
\cite[Thm~D]{choi_lacini_park_sheridan_sings_and_syz_of_secant_vars}.
We expect similar applications and smooth equivalences for the case of arbitrary $k$.
To some extent this is addressed in work in preparation
  \cite{nisiabu_jabu_abcd}, \cite{jagiella_mgr}.
Moreover, Theorem~\ref{thm_cactus_and_socle} motivates  \cite{galazka_keneshlou_sivic_Hilb_of_low_deg} which investigates the loci in the Hilbert scheme parameterising schemes with low socle dimension.

\subsection*{Overview}
After a brief preparation on subschemes of finite schemes in Section~\ref{sec_finite}, we prove Theorem~\ref{thm_cactus_and_socle} in Section~\ref{sec_cacti}.
Finally, Section~\ref{sec_effective} discusses conjectural  optimality of the main result and properties of schemes and algebras with low socle dimensions.

\subsection*{Acknowledgements}
We thank Joachim Jelisiejew for suggesting some references, in particular pointing out an alternative proof of Lemma~\ref{lem_cactus_of_finite_scheme_with_large_socle}.
We also thank the anonymous reviewer for his suggestions and questions.
J.~Buczy{\'n}ski is supported by the project
``Advanced problems in contact manifolds, secant varieties, and their
generalisations (\mbox{APRICOTS+})'', project number
2023/51/B/ST1/02799, awarded by National Science Center, Poland.
M.~Ga{\l}\k{a}zka is supported by the project ``Global Basic Research Laboratory: Algebra and Geometry of Spaces of Tensors, and
Applications'', RS-2024-00414849, awarded by the NRF of Korea.

\section{Subschemes of a finite scheme}
\label{sec_finite}

We prove Theorem~\ref{thm_cactus_and_socle} by reducing it to the special case of $X=R$, where $R$ is finite of degree $r$.
Before that, in a series of lemmas we treat this special case and some related mathematics. In particular, in this section we investigate finite subchemes of a finite scheme, which might be of interest on its own. The statements here are slightly stronger than we need in the course of the proof of main theorem.

If $R$ is a finite scheme, by $\Hilb_{k} R$ we denote the Hilbert scheme of finite subschemes of $R$ of degree $k$.
That is, the $\kk$-points of $\Hilb_{k} R$ correspond uniquely to subschemes $R'\subset R$ of degree $k$, and $\Hilb_{k} R$ is universal with respect to flat deformations of such schemes.
In particular, $\Hilb_{k} R$ comes with the universal subscheme $\ccR'\subset R\times_{\kk} \Hilb_{k} R$ which is flat over $ \Hilb_{k} R$. Moreover, any finite flat family of subschemes of $R$ of degree $k$ is a pullback of $\ccR'\to  \Hilb_{k} R$.

\begin{lem}
  \label{lem_subschemes_of_a_finite_local_scheme}
   Suppose $R=\Spec A$ is a finite local scheme of degree $r$,
     with $A$ a finite local algebra with socle $\gots\subset A$.
   Then $\reduced{(\Hilb_{r-1} R)} = \PP(\gots)$,
   and thus each degree $r-1$ subscheme of $R$ is uniquely determined by a $1$-dimensional $\kk$-linear subspace of $\gots$.
\end{lem}

\begin{prf}
   The bijection of $\kk$-points of  $\reduced{(\Hilb_{r-1} R)}$ and $\PP(\gots)$ is clear: a subscheme $R' \subset R$ corresponds uniquely to an ideal $I_{R'}\subset A$,
   and
   \[
      \deg R'= r-1 \iff \dim I_{R'} = 1 \iff \dim I_{R'} =1 \text{ and } I_{R'} \subset \gots \iff I_{R'} \in \PP(\gots).
   \]
   We now turn this bijection into an isomorphism of reduced schemes by constructing two morphisms --- one each way ---
     that realise this bijection.

   Define $\ccI\subset \ccO_{\PP(\gots)}\otimes_{\kk} A$ to be the universal ideal sheaf: over a $p\in \PP(\gots)$ the ideal is the one dimensional subspace $p\subset \gots$.
   This ideal $\ccI$ defines a family $\ccR'\subset \PP(\gots) \times R$ which is flat, finite and of degree $r-1$ over $\PP(\gots)$.
   Thus we obtain a morphism $\PP(\gots) \to \Hilb_{r-1} R$,
     whose image is contained in $\reduced{(\Hilb_{r-1} R)}$.
   By the above discussion this morphism is bijective
     on $\kk$-points.

   Now we construct a map $\Hilb_{r-1} R \to \PP(A)$.
   Any such map is determined by a line subbundle $L\subset A\otimes_{\kk} \ccO_{\Hilb_{r-1} R}$,
     so that $L$ becomes the pull-back of $\ccO_{\PP (A)}(-1)$ under this morphism.
   Let $L$ be the universal ideal sheaf
     over $\ccO_{\Hilb_{r-1} R}$, that is, if $\ccR\subset R\times_{\kk} \Hilb_{r-1} R$ is the universal family such that the projection $\ccR\to \Hilb_{r-1} R$ is flat finite and of degree $r$, then $L\subset \ccO_R \otimes_{\kk} \ccO_{\Hilb_{r-1} R}$ is its ideal.
   Then $L$ is flat and coherent over $\Hilb_{r-1} R$,
   hence locally free, and moreover its rank is $1$.
   Thus $L$ is the desired line bundle.
   Also this morphism agrees on $\kk$-points with the bijection described in the beginning of the proof.
   Therefore the restriction to $\reduced{(\Hilb_{r-1} R)}$ has its image contained in $\PP(\gots)$ and the morphism is precisely the inverse of $\PP(\gots) \to \reduced{(\Hilb_{r-1} R)}$.
   Hence these two schemes are isomorphic.
\end{prf}

\begin{funfact}
   Although it might be tempting to conjecture that in the setting of Lemma~\ref{lem_subschemes_of_a_finite_local_scheme}
   the Hilbert scheme
   $\Hilb_{r-1} R$
   is also a projective space, this fails already on the simplest nonreduced scheme:
   let $R=\Spec \kk[t]/(t^2)$,
     so that $r=2$, and $\gots=(t)$.
   Then $\Hilb_{1} R= R \neq \PP(\gots) = \Spec \kk$.
   Thus $\Hilb_{r-1} R$ is not necessarily reduced.
   In this case, the map $\Hilb_{1} R\to \PP A\simeq \PP^1$ from the proof of Lemma~\ref{lem_subschemes_of_a_finite_local_scheme}
     is an embedding.
   It would be interesting to see if $\Hilb_{r-1} R\to \PP A$ is an embedding for any finite local $R$.
\end{funfact}

\begin{lem}
  \label{lem_subschemes_of_a_finite_scheme}
   Suppose $R$ is a finite scheme of degree $r$,
    and that $R = \bigsqcup R_i$ with each $R_i$ a finite local scheme of degree $r_i$.
   Then $\Hilb_{r-1} R = \bigsqcup_{i} \Hilb_{r_i-1} R_i$.
\end{lem}

\noprf

\section{Grassmann cacti}
\label{sec_cacti}

In the course of the remaining proofs we denote
\[
 \cactusop{r,k}{X} := \set{E\in \Gr(k, V) \mid
    \exists R\subset X \text{ s.t.} \dim R = 0, \deg R \leqslant r, E \subset \linspan{R} },
\]
so that $\cactus{r,k}{X} = \overline{\cactusop{r,k}{X}}$.

\begin{lem}
  \label{lem_cactus_of_finite_scheme}
   Suppose $R\subset \PP(V)$ is a finite scheme of degree $r$.
   Then the following subsets of $\Gr(k, V)$ are equal:
   \[
      \cactus{r,k}{R} = \cactusop{r,k}{R}
                      = \Gr(k,\linspan{R}).
   \]
\end{lem}

\begin{prf}
   Since $\Gr(k,\linspan{R})$ is closed in $\Gr(k, V)$  it is enough to prove the second equality, which is straightforward,
   since the unique subscheme of $R$ of dimension $0$ and degree $r$ is $R$ itself.
\end{prf}

\begin{lem}
  \label{lem_cactus_of_finite_scheme_with_large_socle}
   Suppose $R\subset \PP(V)$ is a finite scheme of degree $r$
      and that it fails to be of socle dimension at most $k$.
   Then
   \[
       \Gr(k,\linspan{R}) =
      \bigcup_{R'\subsetneq R} \Gr(k,\linspan{R'}),
   \]
   where the union is over all subschemes $R'$ of $R$ of degree $r-1$.
\end{lem}
 A version of this lemma for $R$ local scheme
   (which is the main case)
   but with different proof is also in  \cite[Prop.~4.2]{jelisiejew_keneshlou_k_regular_map_to_Grassmannians}.

\begin{prf}
  Note that $\supset$ part of the claim is automatic.

  To prove the $\subset$ part,
  write $R = \Spec A_0 \sqcup \bigsqcup \Spec A_i$ with each $A_i$ a finite local algebra, and suppose the socle dimension of $A_0$ is $l \geqslant k+1$.
  Let $\gots\subset A_0$ be the socle of $A_0$.
  Define $Q\subset R$ to be the subscheme of degree $r-l$ with
  $Q = \Spec (A_0/\gots) \sqcup \bigsqcup \Spec A_i$.
  We claim that
  \begin{equation}
      \label{equ_Grassmann_inclusion}
      \Gr(k,\linspan{R}) \subset
      \bigcup_{Q\subset R'\subsetneq R} \Gr(k,\linspan{R'}),
  \end{equation}
  so in the union it is enough to use only those $R'$ of degree $r-1$, that contain $Q$.

  If there exists
     $R' \subsetneq R$ as above such that $\linspan{R'} = \linspan{R}$,
  then the inclusion \eqref{equ_Grassmann_inclusion} holds trivially.

  Suppose there exist two different schemes $R'$ and $R''$ both of degree $r-1$ and such that
  $Q\subset R' \subsetneq R$, $Q\subset R'' \subsetneq R$, and  $\linspan{R'} = \linspan{R''}$.
  Then $R$ is the scheme theoretic union $R'\cup R''$, and
  \[
    \linspan{R} = \linspan{R'\cup R''} \subset
    \linspan{R'} + \linspan{R''} = \linspan{R'}.
  \]
  Thus we are back in the first trivial case.

  The final case to consider
    is when for any pair $R', R''\subset R$ of distinct subschemes, both containing $Q$, their linear spans $\linspan{R'}\neq\linspan{R''}$  are different hyperplanes in $\linspan{R}$.
  By Lemmas~\ref{lem_subschemes_of_a_finite_local_scheme} and \ref{lem_subschemes_of_a_finite_scheme},
  the parameter space of subschemes $R'$ such that $\deg R'=r-1$ and $Q\subset R' \subsetneq R$ is equal to
  $\PP(\gots)\simeq \PP^{l-1}$.
  Setting $m=\dim\linspan{R}$,
  we have an injective map
  \[
     \PP(\gots) \to  \Gr(m-1, \linspan{R}) = \PPof{\linspan{R}^*}
  \]
  which sends the scheme $R'$ to its linear span.
  Moreover each element of the image of this injective map
    contains $\linspan{Q}$.
  The Grassmannian of linear subspaces $F\subset \linspan{R}$ such that $\linspan{Q} \subset F$ and $\dim F=m-1$ has dimension at most $l-1=(\deg R - \deg Q)-1$ and it is isomorphic to a projective space, namely the appropriate linear subspace of $\PPof{\linspan{R}^*}$.
  By dimension count, this linear subspace must be equal to the image of the above map
  $\PP(\gots) \to  \PP(\linspan{R}^*)$.
  Therefore, every linear subspace $F$ as above is the linear span of a subscheme $R'$ with $\deg R'=r-1$ and $Q\subset R' \subset R$.

  In order to prove \eqref{equ_Grassmann_inclusion}
  pick any $E \in \Gr(k, \linspan{R})$.
  Define $F'=E+ \linspan{Q}$ and note that $\dim F' \leqslant k + \deg Q \leqslant r-1$ since $l > k$.
  Let $F\subset \linspan{R}$ be any linear subspace such that
    $F' \subset F$ and $\dim F= r-1$.
  We have $\linspan{Q}\subset F$,
    so by the above argument there exists a subscheme $R' \subset R$
    such that $\deg R'=r-1$, $Q\subset R'$ and $\linspan{R'}=F$.
  Thus $E\subset \linspan{R'}$ and therefore $E \in \Gr(k, \linspan{R'})$, as claimed.
  This concludes the proof of the lemma.
\end{prf}

\begin{prf}[ of Theorem~\ref{thm_cactus_and_socle}]
  The inclusion of
  \begin{equation}
    \label{equ_cactus_via_low_socle}
    \overline{\set{E\in \Gr(k, V) \mid
      \exists R\subset X \text{ s.t.} \dim R = 0, \deg R \leqslant r, \socdim R \leqslant k, E \subset \linspan{R}}}
  \end{equation}
  in $\cactus{r,k}{X}$ is automatic.
  To prove the opposite inclusion
  take any $E$ in $\cactusop{r,k}{X}$,
    so that there exists a subscheme $R \subset X$ with $\dim R = 0$, $\deg R \leqslant r$, $E \subset \linspan{R}$.
  Suppose without loss of generality that $R$ is a minimal such scheme,
    that is for any $R'\subset R$, if $E \subset \linspan{R'}$ then $R'=R$.
  If $R$ fails to have socle dimension at most $k$,
  then
  by Lemmas~\ref{lem_cactus_of_finite_scheme} and \ref{lem_cactus_of_finite_scheme_with_large_socle}
   there exists a smaller scheme $R'\subset R$ with $\deg R'< \deg R$ and
   $E \subset \linspan{R'}$, a contradiction with minimality of $R$.
  Thus $R$ has socle dimension at most $k$,
    and therefore $E$ is in the set~\eqref{equ_cactus_via_low_socle}.
  Therefore $\cactusop{r,k}{X}$ is contained in the set~\eqref{equ_cactus_via_low_socle}, which is closed,
  thus also $\cactus{r,k}{X}$ is contained in the same set, concluding the proof of the theorem.
\end{prf}

\section{Effectiveness}
\label{sec_effective}

In this section we briefly discuss additional questions posed by the reviewer of this article. We do not have precise answers yet, but the topics are being independently investigated and lay outside the scope of this article.

Theorem~\ref{thm_cactus_and_socle} states that to define
$\cactus{r,k}{X}$ it is enough to consider only finite schemes of socle dimension at most $k$.
But is it optimal?
Perhaps it is enough to use even fewer schemes? For instance, are Gorenstein schemes enough to define $\cactus{r,2}{X}$?
In general, the answer has not yet been written down in the literature.
\begin{conj}
   Suppose $X$ is a projective variety and $r,k$ are positive integers. If $X\subset \PP V$ is a ``sufficiently ample'' embedding, then the irreducible components of $\cactus{r,k}{X}$ are in a bijection with irreducible components of the Hilbert scheme $\Hilb_r(X)$ whose general element has socle dimension at most $k$.
\end{conj}
The notion of sufficiently ample embedding is discussed in detail in \cite{sidman_smith_determinantal}, \cite{nisiabu_jabu_farnik_cactus_Fujita}, and references therein.
In particular, if $X\subset \PP^n$ is any embedding, then the $d$-th Veronese reembedding $\nu_d \colon X\hookrightarrow \PP^N$ for $d\gg 0$ is sufficiently ample in this sense.

If the conjecture is true, then Theorem~\ref{thm_cactus_and_socle} is optimal, and one cannot remove any further large families of schemes in general.
The conjecture is known to experts in the case $k=1$: it follows from \cite[Thm~1.6]{nisiabu_jabu_cactus} and the methods in \cite{nisiabu_jabu_farnik_cactus_Fujita}.
More explicitly,
  it is stated for $X=\PP^n$
  in \cite[Thm~1.2]{jabu_keneshlou_cactus_scheme}.
For larger values of $k>1$, the only case known in the literature is $\cactus{8,3}{\nu_d(\PP^n)}$ for $d \geqslant 5$, see \cite[\S6]{galazka_mandziuk_rupniewski_distinguishing}.

Another open problem is to describe the loci in the Hilbert scheme $\Hilb_r^{\socdim\leqslant k}(X)$ consisting of schemes that have socle dimension at most $k$.
This locus is open in the Hilbert scheme
(similarly to the openness of $\Hilb_r^{\Gor}(X)=\Hilb_r^{\socdim\leqslant 1}(X)$).
In particular, it is not known what is the smallest $r$, for which $\Hilb_r^{\socdim\leqslant 2}(X)$ is reducible.
If $\cchar \kk=0$, based on \cite{galazka_keneshlou_sivic_Hilb_of_low_deg} we expect it should be $r=13$.

More generally, the properties and structure of finite local Gorenstein algebras are relatively well understood, see for instance \cite[\S2.1,2.3, 3.3, 3.4]{jelisiejew_PhD} and references therein.
The famous Gorenstein symmetry of (pieces of) Hilbert functions or self-duality increase the research interest in these algebras.
On the other hand, finite local algebras that have socle dimension $k$ for $k\geqslant 2$ are less understood.
Important examples include level algebras, see for instance
\cite{boij_Betti_numbers_of_compressed_level_algebras},
\cite{cho_iarrobino_Hilbert_functions_and_level_algebras},
\cite{chipalkatti_geramita_parameter_spaces_for_Artin_level_algebras},
\cite{geramita_harima_migliore_shin_Hilbert_function_of_level_algebra}.

\addcontentsline{toc}{section}{References}

\bibliography{grassmann_cactus_and_socle.bbl}
\end{document}